\documentclass[11pt]{article}

\usepackage{amsmath}
\usepackage{amssymb}
\usepackage{amsbsy}
\usepackage{amsfonts}
\usepackage{mathtools}
\usepackage{bbm}
\usepackage{comment}
\newtheorem{theorem}{Theorem}[section]%
\newtheorem{lemma}[theorem]{Lemma}%
\newenvironment{pf}{\medskip\noindent{Proof:}
  \hspace{-.5cm}      \enspace}{\hfill \qed \newline \smallskip}

\def\det{{\rm det}}

\def\f{\noindent}

\def\mod{\hbox{\rm mod}\;}

\newcommand{\qed}{\mbox{\raisebox{0.7ex}{\fbox{}}} \vspace{4truemm}}

\begin{document}

\begin{center}
{\Large{\textbf{Complexity of the circulant foliation over a graph. }}}
\end{center}

\vskip 5mm {\small
\begin{center}
{\textbf{Y.~S.~Kwon,}}\footnote{{\small\em Department of
Mathematics, Yeungnam University, Korea}}
{\textbf{A.~D.~Mednykh,}}\footnote{{\small\em Sobolev Institute of
Mathematics, Novosibirsk State University, Russia}} {\textbf{I.~A.~Mednykh,}}\footnote{{\small\em Sobolev
Institute of Mathematics, Novosibirsk State University, Russia}}
 \end{center}}

\title{ \vspace{-1.2cm}
Complexity of the  circulant foliation over a graph.
\thanks{Supported by }}

\begin{abstract}
In the present paper, we investigate the complexity of infinite family of graphs $H_n=H_n(G_1,\,G_2,\ldots,G_m)$ obtained as a circulant foliation over a graph $H$ on $m$ vertices with fibers $G_1,\,G_2,\ldots,G_m.$ Each fiber $G_i=C_n(s_{i,1},\,s_{i,2},\ldots,s_{i,k_i})$ of this foliation is the circulant graph on $n$ vertices with jumps $s_{i,1},\,s_{i,2},\ldots,s_{i,k_i}.$ This family includes the family of generalized Petersen graphs, $I$-graphs, sandwiches of circulant graphs, discrete torus graphs and others.

We obtain a closed formula for the number $\tau(n)$ of spanning trees in $H_n$ in terms of Chebyshev polynomials, investigate  some arithmetical
  properties of this function  and find its asymptotics  as $n\to\infty.$
\bigskip

\f\textbf{AMS classification:} 05C30, 39A10\\ \textbf{Keywords:}
spanning tree,  circulant graph, Petersen graph, $I$-graph, Laplacian matrix,
Chebyshev polynomial
\end{abstract}

\section{Introduction}
 Let $G$ be a finite connected graph. By the {\it complexity} $\tau(G)$ of the graph $G$  we mean the number of its spanning trees.  
 The complexity is very important algebraic invariant of a graph. 
 Various approaches to its computation are given in the papers \cite{BoePro, Golin2010, SuWanZhang, XiebinLinZhang, AD11, CCY07, KwonMedMed}. 
 For an infinite family of graphs $G_n,\, n\in \mathbb{N}$ one can introduce complexity function $\tau(n)=\tau(G_n).$ 
  In statistical physics \cite{Wu77, SW00, GutRog, Louis}, it is important to know the behavior of the function $\tau(n)$ for sufficiently large values of $n$.
 
The aim of the present paper is to investigate analytical, arithmetical  and  asymptotic properties of complexity function for circulant foliation over a given graph.  We note that this family is quite rich. It includes circulant graphs, generalized Petersen graphs, $I$-, $Y$-, $H$-graphs, discrete tori and others.

The structure of the paper is as follows. Some preliminary results and basic definitions are given in Section~\ref{basic}. In  Section~\ref{circulant foliation} we define the notion of circulant foliation over a graph. In Section~\ref{counting},
we present explicit formulas for the number of spanning trees of graphs $H_n=H_n(G_1,\,G_2,\ldots,G_m)$ obtained as a circulant foliation over a graph $H$ on $m$ vertices with fibers $G_1,\,G_2,\ldots,G_m.$ Each fiber $G_i=C_n(s_{i,1},\,s_{i,2},\ldots,s_{i,k_i})$ of this foliation is the circulant graph on $n$ vertices with jumps $s_{i,1},\,s_{i,2},\ldots,s_{i,k_i}$. The formulas will be given in terms of Chebyshev polynomials.  In Section~\ref{arithmetic},
we provide some arithmetical properties of the complexity function for the family $H_n.$ More precisely, we show
that the number of spanning trees in  the graph $H_n$ can be represented in
the form $\tau(n)=p \,n\,\tau(H) \,a(n)^2,$ where $a(n)$ is an integer sequence and $p$ is a prescribed
natural number depending   on jumps and the parity of  $n.$ In Section~\ref{asymptotic}, we use
explicit formulas for the complexity in order to produce its asymptotic. In the last section,
we illustrate the obtained results by a series of examples.

\section{Basic definitions and preliminary facts}\label{basic}

Consider a connected finite graph $G,$ allowed to have multiple edges but without loops. We denote the vertex and edge set of $G$ by $V(G)$ and $E(G),$ respectively.
 Given $u, v\in V(G),$ we set $a_{uv}$ to be the number of edges between vertices $u$ and $v.$ 
 The matrix $A=A(G)=\{a_{uv}\}_{u, v\in V(G)}$ is called \textit{the adjacency matrix} of the graph $G.$ 
 The degree $d_v$ of a vertex $v \in V(G)$ is defined by $d_v=\sum_{u\in V(G)}a_{uv}.$ Let $D=D(G)$ be the diagonal matrix indexed by the elements of $V(G)$ with $d_{vv} = d_v.$ 
 The matrix $L=L(G)=D(G)-A(G)$ is called \textit{the Laplacian matrix}, or simply \textit{Laplacian}, of the graph $G.$  
 Let  $X=\{x_v,\,v\in V(G)\}$ be the set of variables and let $X(G)$ be the diagonal matrix indexed by the elements of $V(G)$ with diagonal elements $x_v.$ 
 Then the {\it generalized Laplacian matrix} of $G,$ denoted by $L(G,X),$ is given by $L(G,X)=X(G)-A(G).$ In the particular case $x_v=d_v,$  we have $L(G,X)=L(G).$ 
In what follows, by $I_n$ we denote the identity matrix of order $n.$

We call an $n\times n$ matrix {\it circulant,} and denote it by $circ(a_0, a_1,\ldots,a_{n-1})$ if it is of the form
$$circ(a_0, a_1,\ldots, a_{n-1})=
\left(\begin{array}{ccccc}
a_0 & a_1 & a_2 & \ldots & a_{n-1} \\
a_{n-1} & a_0 & a_1 & \ldots & a_{n-2} \\
  & \vdots &   & \ddots & \vdots \\
a_1 & a_2 & a_3 & \ldots & a_0\\
\end{array}\right).$$

Recall \cite{PJDav} that the eigenvalues of matrix $C=circ(a_0,a_1,\ldots,a_{n-1})$ are given by the following simple formulas $\lambda_j=p(\varepsilon^j_n),\,j=0,1,\ldots,n-1$
 where $p(x)=a_0+a_1 x+\ldots+a_{n-1}x^{n-1}$ and $\varepsilon_n$ is an order $n$ primitive root of the unity. Moreover, the circulant matrix $C=p(T_n ),$ 
 where $T_{n}=circ(0,1,0,\ldots,0)$ is the matrix representation of the shift operator $T_{n}:(x_0,x_1,\ldots,x_{n-2},x_{n-1})\rightarrow(x_1, x_2,\ldots,x_{n-1},x_0).$ 
 For any $i=0,\ldots, n-1$, let ${{\bf v}_i} =(1,\varepsilon_n^i ,\varepsilon_n^{2i},\ldots,\varepsilon_n^{(n-1)i})^t$ be a column vector of length $n.$ 
 We note that all $n\times n$ circulant matrices share the same set of linearly independent eigenvectors $\textbf{v}_0, \textbf{v}_1, \ldots, \textbf{v}_{n-1}.$
Hence, any set of $n\times n$ circulant matrices can be simultaneously diagonalizable.

Let $s_1,s_2,\ldots,s_k$ be integers such that $1\leq s_1<s_2<\ldots<s_k\leq\frac{n}{2}.$ 
The graph $C_n(s_1,s_2,\dots,s_k)$ with $n$ vertices $0,1,2,\dots,n-1$ is called {\it circulant graph} if the vertex $i,0\leq i\leq  n-1$ is adjacent to the vertices $i\pm s_1, i\pm s_2,\dots,i\pm s_k$ (mod $n$). 
All vertices of the graph are of even degree $2k$. If $n$ is even and $s_k=\frac{n}{2}$, then the vertices $i$ and $i+s_k$ are connected by two edges. 
In this paper, we also allow the empty circulant graph $C_n(\emptyset)$ consisting of $n$ isolated vertices.

\section{Circulant foliation over a graph}\label{circulant foliation}

Let $H$ be a connected finite graph on vertices $v_1,v_2,\ldots,v_m$, allowed to have multiple edges but without loops. 
Denote by $a_{i\,j}$ the number of edges between vertices $v_i$ and $v_j.$ Since $H$ has no loops, we have $a_{i\,i}=0.$ 
To define the circulant foliation $H_n=H_n(G_1,\,G_2,\ldots,G_m)$ we prescribe  to each vertex $v_i$ a circulant graph $G_i=C_n(s_{i,1},\,s_{i,2},\ldots,s_{i,k_i}).$ 
Then the {\it circulant foliation} $H_n=H_n(G_1,\,G_2,\ldots,G_m)$ over $H$ with fibers $G_1,\,G_2,\ldots,G_m$ is a graph with the vertex set $V(H_n)=\{(k,\,v_i)\ | \,k=1,2,\ldots n,\,i=1,2,\ldots,m\},$ 
where for a fixed $k$ the vertices $(k,\,v_i)$ and $(k,\,v_j)$ are connected by $a_{i\,j}$ edges, while for a fixed $i,$ the vertices $(k,\,v_i),\,k=1,2,\ldots n$ 
form a graph $C_n(s_{i,1},\,s_{i,2},\ldots,s_{i,k_i})$ in which the vertex $(k,\,v_i)$ is adjacent to the vertices $(k\pm s_{i,1},v_i),(k\pm s_{i,2},\,v_i),\ldots,(k\pm s_{i,k_i},\,v_i),(\textrm{mod}\ n).$

There is a projection $\varphi: H_n\to H$ sending  the vertices $(k,\,v_i),\,k=1,\ldots,n$ and edges between them to the vertex $v_i$ and for given $k$ each edge between the vertices $(k,\,v_i)$ 
and $(k,\,v_j),\,i\neq j$ bijectively to an edge between $v_i$ and $v_j.$ For each vertex $v_i$ of graph $H$ we have $\varphi^{-1}(v_i)=G_i,\,i=1,2,\ldots,m.$ 
 
Consider an action of the cyclic group $\mathbb{Z}_n$ on the graph $H_n$ defined by the rule $(k,\,v_{i})\to(k+1,\,v_{i}),\,k\,\mod n.$ 
Then the group $\mathbb{Z}_n$ acts free on the set of vertices and the set of edges and the factor graph $H_{n}/\mathbb{Z}_{n}$ is an {\it equipped  graph} $\widehat{H}$ obtained from the graph $H$ 
by attaching $k_i$ loops to each $i$-th vertex of $H.$

By making use of the voltage technique \cite{GrossTucker}, one  can construct the graph $H_n$ in the following way. 
We put an orientation to all edges of $\widehat{H}$ including loops. Then we prescribe the voltage $0$ to all edges of subgraph $H$ of $\widehat{H}$ and the voltage $s_{i,j},\,\mod n$ to 
the $j$-th loop attached to $i$-th vertex of $H.$ The respective voltage covering is the graph $H_n.$ It is well known  that the obtained graph $H_n$ is connected if and only if 
the voltages $\{s_{i,j},\,\mod n\}$ generate the full group $\mathbb{Z}_n.$ Equivalently, $H_n$ is connected if and only if $\gcd(n,s_{i,j},\,i=1,\ldots,m,\,j=1,\ldots, k_i)=1.$ 
Moreover, if $r$ is a unit in the ring $\mathbb{Z}_n$ (that is, there is an element $r^{\prime}$ in $\mathbb{Z}_n$ such that $r r^{\prime}=1$ ), then the graphs $H_n$ and $H_n^{\prime}$ 
obtained by the voltage assignments $\{s_{i,j},\,\mod n\}$ and $\{r\,s_{i,j},\,\mod n\}$ are isomorphic.

Recall that the adjacency matrix of the circulant graph $C_n(s_1,s_2,\ldots,s_k)$ on the vertices $1,2,\ldots,n$ has the form $\sum\limits_{p=1}^{k}(T_{n}^{s_p}+T_{n}^{-s_p}).$ 
Let the adjacency matrix of the graph $H$ be
$$A(H)=
\left(\begin{array}{ccccc}
0 & a_{1,2} & a_{1,3} & \ldots & a_{1,m} \\
a_{2,1} & 0 & a_{2,3} & \ldots & a_{2,m} \\
  & \vdots &   & \ddots & \vdots \\
a_{m,1} & a_{m,2} & a_{m,3} & \ldots & 0\\
\end{array}\right).$$
Then, the adjacency matrix of the circulant foliation $H_n=H_n(G_1,\,G_2,\ldots,G_m)$ over a graph $H$ with fibers $G_i=C_n(s_{i,1},\,s_{i,2},\ldots,s_{i,k_i}),\,i=1,2,\ldots,n$ is given by
$$A(H_n)=
\left(\begin{array}{ccccc}
\sum\limits_{p=1}^{k_1}(T_{n}^{s_{1,p}}+T_{n}^{-s_{1,p}}) & a_{1,2}I_n & a_{1,3}I_n & \ldots & a_{1,m}I_n \\
a_{2,1}I_n & \sum\limits_{p=1}^{k_2}(T_{n}^{s_{2,p}}+T_{n}^{-s_{2,p}}) & a_{2,3}I_n & \ldots & a_{2,m}I_n \\
  & \vdots &   & \ddots & \vdots \\
a_{m,1}I_n & a_{m,2}I_n & a_{m,3}I_n & \ldots & \sum\limits_{p=1}^{k_m}(T_{n}^{s_{m,p}}+T_{n}^{-s_{m,p}})\\
\end{array}\right).$$

\bigskip
As the first example,   we consider the {\it sandwich graph} $SW_n=H_n(G_1,\,G_2,\ldots,G_m)$ formed by the circulant graphs $G_i=C_n(s_{i,1},\,s_{i,2},\ldots,s_{i,k_i}).$ 
To create $SW_n$ we take $H$ to be the path graph on $m$ vertices $v_1,v_2, \ldots,v_m$ with the end points $v_1$ and $v_m.$  
 A very particular case of this construction, known as $I$-graph $I(n,k,l),$ occurs  by taking $m=2,\,G_1=C_n(k)$ and $G_2=C_n(l).$ 
 Also, the generalized Petersen graph \cite{SS09} arises as $GP(n,k)=I(n,k,1).$ The sandwich  of two circulant graphs $H_n(G_1,G_2)$ was investigated in \cite{AbrBaiMed}.

As the second example, we consider the {\it generalized $Y$-graph} $Y_n=Y_n(G_1,G_2,G_3),$ where $G_1,G_2,G_3$ are given circulant graphs on $n$ vertices. 
To construct $Y_n,$ we consider a  $Y$-shape graph $H$ consisting of four vertices $v_1,v_2,v_3,v_4$ and three edges $v_1v_4,v_2v_4,v_3v_4.$ Let $G_4=C_n(\emptyset)$ be the empty graph  of $n$ on vertices. Then, by definition, we put $Y_n=H_n(G_1,G_2,G_3,G_4).$ In a particular case, $G_1=C_n(k),G_2=C_n(l),$ and $G_3=C_n(m),$ the graph $Y_n$ coincides with the $Y$-graph $Y(n;k,l,m)$ defined earlier in \cite{BigIYH, HorBou}.  
 
The third example is the  {\it generalized $H$-graph} $H_n(G_1,G_2,G_3,G_4,G_5,G_6),$ where $G_1,G_2,G_3,G_4$ are given circulant graphs  and $G_5= G_6=C_n(\emptyset)$ are the empty graphs on $n$  vertices. 
In this case, we take $H$ to be the graph with vertices $v_1,v_2,v_3,v_4,v_5,v_6$ and edges $v_1v_5,v_5v_3,v_2v_6,v_6v_4,v_5v_6.$ 
In the case $G_1=C_n(i),G_2=C_n(j),G_3=C_n(k),G_4=C_n(l),$ we get the graph $H(n;i,j,k,l)$ investigated in the paper \cite{HorBou}.
Shortly, we will write $H_n(G_1,G_2,G_3,G_4)$  ignoring the last two empty graph entries.

\section{Counting the number of spanning trees in the graph $H_n$}\label{counting}

Let $H$ be a finite connected graph with the vertex set $V(H)=\{v_1,\,v_2,\ldots,v_m\}.$ Consider the circulant foliation $H_n=H_n(G_1,\,G_2,\ldots,G_m),$ 
where $G_i=C_n(s_{i,1},\,s_{i,2},\ldots,s_{i,k_i}),\,i=1,2,\ldots,m.$ Let $L(H,\,X)$ be the generalized Laplacian of graph $H$ with the set of variables $X=(x_1,x_2,\ldots,x_m).$ 
We specify $X$ by setting $x_{i}=2k_{i}+d_i-\sum\limits_{p=1}^{k_i}(z^{s_{i,p}}+z^{-s_{i,p}})$ and put $P(z)=\det(L(H,\,X))$, where $d_i$ is the degree of $v_i$ in $H$.  
We note that $P(z)$ is an integer Laurent polynomial.  
Consider one more specification $L(H,\,W)$ for generalized Laplacian of $H$  with the set  $W=(w_1,w_2,\ldots,w_m),$ where $w_{i}=2k_{i}+d_i-\sum\limits_{p=1}^{k_i}2\,T_{s_{i,p}}(w)$ 
and $T_k(w)=\cos(k\arccos w)$ is the Chebyshev polynomial of the first kind. See \cite{MasHand} for the basic properties of the Chebyshev polynomials. 
We set $Q(w)=\det(L(H,\,W)).$ Then $Q(w)$ is an integer polynomial of degree $s=s_{1,k_1}+s_{2,k_2}\ldots+s_{m,k_m}.$ For our convenience, we will call $Q(w)$ a {\it Chebyshev transform}  of $P(z).$ 
The following lemma holds.

\bigskip
\begin{lemma}\label{lemma1}
We have $P(z)=Q(w)$ with $w=\frac{1}{2}(z+\frac{1}{z})$ and $Q(w)$ is the order $s$ polynomial with the leading coefficient $(-1)^{m}2^s,$ where $s=\sum\limits_{i=1}^{m}s_{i,k_i}.$ 
Moreover, $$Q(1)=0,\,Q^{\prime}(1)=-2q\tau(H)\neq 0^{},$$ where $q=\sum\limits_{i=1}^{m}\sum\limits_{j=1}^{k_{i}}s_{i,j}^2$ and $\tau(H)$ is the number of spanning trees in the graph $H.$
In particular, $Q(w)$ has a simple  root $w=1$  and $P(z)$ has a double  root $z=1.$
\end{lemma}

\begin{pf}
The equality $P(z)=Q(w)$ follows from the identity $T_{n}(\frac{1}{2}(z+\frac{1}{z}))=\frac{1}{2}(z^{n}+\frac{1}{z^{n}}).$ Recall that the leading term of $T_{n}(w)$ is $2^{n-1}w^{n}.$ The leading term of $Q(w)$ is coming from the product $\prod\limits_{i=1}^{m}(-2T_{s_{i,k_{i}}}(w))$ and is equal to $(-1)^{m}2^{s}w^{s},$ where $s=\sum\limits_{i=1}^{m}s_{i,k_{i}}.$

Let $a_{i,j}$ be the number of edges between $i$-th and $j$-th vertices of the graph $H.$ Then \begin{equation}\label{formula1}Q(w)=\det\left(\begin{array}{ccccc}x_1 & -a_{1,2} & -a_{1,3} & \ldots & -a_{1,m} \\
-a_{2,1} & x_2 & -a_{2,3} & \ldots & -a_{2,m} \\
  & \vdots &   & \ddots & \vdots \\
-a_{m,1} & -a_{m,2} & -a_{m,3} & \ldots & x_m\end{array}\right),\end{equation} where $x_i=x_i(w)=2k_i+d_i-\sum\limits_{j=1}^{k_i}2T_{s_{i,j}}(w),\,i=1,2,\ldots,m.$ In particular, for $w=1$ we have $x_i=d_i.$ Hence, $Q(1)=0$ because of valency of $i$-th vertex is $d_i=\sum_j a_{i,j}.$ Let $x_i^{\prime}=x_i^{\prime}(w)$ be the derivative of $x_i$ with respect to $w.$ Then $$Q^{\prime}(w)=\det\left(\begin{array}{ccccc}
x^{\prime}_1 & -a_{1,2} & -a_{1,3} & \ldots & -a_{1,m} \\
0 & x_2 & -a_{2,3} & \ldots & -a_{2,m} \\
  & \vdots &   & \ddots & \vdots \\
0 & -a_{m,2} & -a_{m,3} & \ldots & x_m\end{array}\right)+ \det\left(\begin{array}{ccccc}
x_1 & 0 & -a_{1,3} & \ldots & -a_{1,m} \\
-a_{2,1} & x^{\prime}_2 & -a_{2,3} & \ldots & -a_{2,m} \\
  & \vdots &   & \ddots & \vdots \\
-a_{m,1} & 0 & -a_{m,3} & \ldots & x_m\end{array}\right)$${}
$$+\ldots+ \det\left(\begin{array}{ccccc}
x_1 & -a_{1,2}& -a_{1,3} & \ldots & 0 \\
-a_{2,1} & x_2 & -a_{2,3} & \ldots & 0 \\
  & \vdots &   & \ddots & \vdots \\
-a_{m,1} & -a_{m,2} & -a_{m,3} & \ldots & x^{\prime}_m\end{array}\right)$${}
$$=x^{\prime}_1(w)Q_{1,1}(w)+x^{\prime}_2(w)Q_{2,2}(w)+\ldots+x^{\prime}_m(w) Q_{m,m}(w),$$
where $Q_{i,i}(w)$ is the $(i,i)$-th minor of the matrix in formula (\ref{formula1}). For $w=1$ this matrix coincides with the Laplacian of $H.$ By the Kirchhoff theorem we have $$Q_{1,1}(1)=Q_{2,2}(1)=\ldots=Q_{m,m}(1)=\tau(H),$$ where $\tau(H)$ is the number of spanning trees in the graph $H.$

Since $T_s^{\prime}(w)=s\, U_{s}(w),$ where $U_s(w)$ is the Chebyshev polynomial of the second kind and $U_s(1)=s,$ we have $x^{\prime}_i(w)=-\sum\limits_{j=1}^{k_i}2s_{i,j}U_{s_{i,j}}(w)$ and $x^{\prime}_i(1)=-2\sum\limits_{j=1}^{k_i}s_{i,j}^2.$

As a result, $Q^{\prime}(1)=(x_1^{\prime}(1)+\ldots+x_m^{\prime}(1))\tau(H)= -2\sum\limits_{i=1}^{m}\sum\limits_{j=1}^{k_i}s_{i,j}^2\tau(H)=-2q\,\tau(H).$ \end{pf}

The main result of this section is the following theorem.
\bigskip

\begin{theorem}\label{theorem1}
The number of spanning trees $\tau(n)$ in the graph $H_{n}(G_1,\,G_2,\ldots,G_m)$ is given by the formula $$\tau(n)=\frac{n\,\tau(H)}{q}\prod_{p=1}^{s-1}|2T_n(w_p)-2|,$$ where $s=s_{1,k_1}+s_{2,k_2}\ldots+s_{m,k_m},\,w_p\,(p=1,2,\ldots,s-1)$ are all the roots different from $1$ of the equation $Q(w)=0$, $\tau(H)$ is the number of spanning trees in the graph $H$ and $q=\sum\limits_{i=1}^{m}\sum\limits_{j=1}^{k_i}s_{i,j}^2.$
\end{theorem}

\begin{pf}
By the classical Kirchhoff theorem, the number of spanning trees $\tau(n)$ is equal to the product of nonzero eigenvalues of the Laplacian of a graph $H_{n}(G_1,\,G_2,\ldots,G_m)$ divided by the number of its vertices $m\times n.$ To investigate the spectrum of Laplacian matrix, we consider the shift operator $T_{n}=circ(0,1,\ldots,0).$ The Laplacian $L=L(H_{n}(G_1,\,G_2,\ldots,G_m))$ is given by the matrix
$$L=
\left(\begin{array}{ccccc}
A_{1}(T_{n}) & -a_{1,2}I_n & -a_{1,3}I_n & \ldots & -a_{1,m}I_n \\
-a_{2,1}I_n & A_{2}(T_{n}) & -a_{2,3}I_n & \ldots & -a_{2,m}I_n \\
 & \vdots & & \ddots & \vdots \\
-a_{m,1}I_n & -a_{m,2}I_n & -a_{m,3}I_n & \ldots & A_{m}(T_{n})\\
\end{array}\right),$$
where $A_{i}(z)=2k_{i}+d_{i}-\sum\limits_{j=1}^{k_{i}}(z^{s_{i,j}}+z^{-s_{i,j}}),\,i=1,\ldots,m.$

The eigenvalues of circulant matrix $T_{n}$ are $\varepsilon_{n}^{j},\,j=0,1,\ldots,n-1,$ where $\varepsilon_n=e^\frac{2\pi i}{n}.$ Since all of them are distinct, the matrix $T_{n}$ is conjugate to the diagonal matrix $\mathbb{T}_{n}=diag(1,\varepsilon_{n},\ldots,\varepsilon_{n}^{n-1})$ with diagonal entries $1,\varepsilon_{n},\ldots,\varepsilon_{n}^{n-1}$. To find spectrum of $L,$ without loss of generality, one can assume that $T_{n}=\mathbb{T}_{n}.$ Then all $n\times n$ blocks of $L$ are diagonal matrices. This essentially simplifies the problem of finding eigenvalues of the block matrix $L.$ Indeed, let $\lambda$ be an eigenvalue of $L$ and let $(x_{1},x_{2},\ldots,x_{m})$ with $x_{i}=(x_{i,1},x_{i,2}\ldots,x_{i,n})^t,\,i=1,\ldots,m$ be the respective eigenvector. Then we have the following system of equations
\begin{equation}\label{equationL1}\left(\begin{array}{ccccc}
A_{1}(\mathbb{T}_{n}) -\lambda\,I_{n}& -a_{1,2}I_{n} & -a_{1,3}I_{n} & \ldots & -a_{1,m}I_{n} \\
-a_{2,1}I_{n} & A_{2}(\mathbb{T}_{n})-\lambda\,I_{n} & -a_{2,3}I_{n} & \ldots & -a_{2,m}I_{n} \\
 & \vdots & & \ddots & \vdots \\
-a_{m,1}I_{n} & -a_{m,2}I_{n} & -a_{m,3}I_{n} & \ldots & A_{m}(\mathbb{T}_{n})-\lambda\,I_{n}\\
\end{array}\right)
\left(\begin{array}{c}
x_{1}\\
x_{2}\\
\vdots \\
x_{m}\\
\end{array}\right)=0.\end{equation}
Recall that all blocks in the matrix under consideration are diagonal
$n\times n$-matrices and the $(j,j)$-th entry of $\mathbb{T}_{n}$ is equal to $\varepsilon_n^{j-1}.$

Hence, the equation (\ref{equationL1}) splits into $n$  equations

\smallskip

\begin{equation}\label{equationL2}\left(\begin{array}{ccccc}
A_{1}(\varepsilon_n^{j}) -\lambda& -a_{1,2}  & -a_{1,3}  & \ldots & -a_{1,m}  \\
-a_{2,1} & A_{2}(\varepsilon_n^{j})-\lambda  & -a_{2,3}  & \ldots & -a_{2,m}  \\
 & \vdots & & \ddots & \vdots \\
-a_{m,1}  & -a_{m,2}  & -a_{m,3}  & \ldots & A_{m}(\varepsilon_n^{j})-\lambda \\
\end{array}\right)
\left(\begin{array}{c}
x_{1,j+1}\\
x_{2,j+1}\\
\vdots \\
x_{m,j+1}\\
\end{array}\right)=0,\end{equation}
$j=0,1,\ldots,n-1$. Each equation gives $m$ eigenvalues of $L,$  say $\lambda_{1,j},\lambda_{2,j},\ldots,\lambda_{m,j}.$ To find these eigenvalues  we set

\smallskip

\begin{equation}\label{equationL23}P(z,\lambda)=\det \left(\begin{array}{ccccc}
A_{1}(z) -\lambda & -a_{1,2}  & -a_{1,3}  & \ldots & -a_{1,m}  \\
-a_{2,1} & A_{2}(z)-\lambda  & -a_{2,3}  & \ldots & -a_{2,m}  \\
 & \vdots & & \ddots & \vdots \\
-a_{m,1}  & -a_{m,2}  & -a_{m,3}  & \ldots & A_{m}(z)-\lambda \\
\end{array}\right)
.\end{equation}

Then $\lambda_{1,j},\lambda_{1,j},\ldots,\lambda_{m,j}$ are  roots of the equation
\begin{equation}P(\varepsilon_n^{j},\lambda)=0.\end{equation}
In particular, by Vieta's theorem, the product $p_{j}=\lambda_{1,j}\lambda_{2,j}\ldots\lambda_{m,j}$ is given by the formula $p_{j}=P(\varepsilon_n^j,0)=P(\varepsilon_n^j),$  where $P(z)$  is the same as in Lemma~\ref{lemma1}.

Now, for any $j=0,\ldots, n-1,$ matrix $L$ has $m$ eigenvalues $\lambda_{1,j},\lambda_{2,j},\ldots,\lambda_{m,j}$ satisfying the order $m$ algebraic equation $P(\varepsilon_n^{j},\lambda)=0.$ In particular, for $j=0$ and $\lambda=\lambda_{i,0},\,i=1,2,\ldots,m$ we have $P(1,\lambda)=0.$ In this case, $A_{i}(1)=d_{i},\,i=1,2,\ldots,m.$ One can see that the polynomial $P(1,\lambda)$ is the characteristic polynomial for Laplace matrix of the graph $H$ and its roots are eigenvalues of $H.$

Note that
 $\lambda_{1,0}=0$ and the product of nonzero eigenvalues $\lambda_{2,0}\lambda_{3,0}\ldots\lambda_{m,0}$ is equal to $m\,\tau(H),$ where $\tau(H)$ is the number of spanning trees in the graph $H.$

Now we have
\begin{equation}\label{tauH}
\tau(n)=\frac{1}{m\,n}\lambda_{2,0}\lambda_{3,0}\ldots\lambda_{m,0} \prod\limits_{j=1}^{n-1}\lambda_{1,j}\lambda_{2,j}\ldots\lambda_{m,j}=
\frac{m\,\tau(H)}{m\,n}\prod\limits_{j=1}^{n-1}p_j=\frac{\tau(H)}{n}\prod\limits_{j=1}^{n-1}P(\varepsilon_n^j).
\end{equation}

To continue the proof we replace the Laurent polynomial $P(z)$ by $\widetilde{P}(z)=(-1)^{m}z^{s}P(z).$
Then $\widetilde{P}(z)$ is a monic polynomial of the degree $2s$ with the same roots as $P(z).$ We note that
\begin{equation}\label{newP}
\prod\limits_{j=1}^{n-1}\widetilde{P}(\varepsilon_{n}^{j})=(-1)^{m(n-1)}\varepsilon_{n}^{\frac{(n-1)n}{2}s} \prod\limits_{j=1}^{n-1}P(\varepsilon_{n}^{j})=(-1)^{(m+s)(n-1)}\prod\limits_{j=1}^{n-1}P(\varepsilon_{n}^{j}).
\end{equation}

By Lemma~\ref{lemma1}, all roots of polynomials $\widetilde{P}(z)$ and $Q(w)$ are $1,1,z_{1},1/z_{1},\ldots,z_{s-1},1/z_{s-1},\,z_{j}\neq1\textrm{ and }1\neq w_{j}=\frac{1}{2}(z_{j}+z_{j}^{-1}),\,j=1,\ldots,s-1,$ respectively. Also, we can recognize the complex numbers $\varepsilon_{n}^{j},\,j=1,\ldots,n-1$ as the roots of polynomial $\frac{z^n-1}{z-1}.$
By the basic properties of resultant (\cite{PrasPol}, Ch.~1.3) we have
\begin{eqnarray}\label{Hlemma}
\nonumber &&\prod\limits_{j=1}^{n-1}\widetilde{P}(\varepsilon_{n}^{j})=\textrm{Res}(\widetilde{P}(z),\frac{z^{n}-1}{z-1}) =\textrm{Res}(\frac{z^{n}-1}{z-1},\widetilde{P}(z))\\
&&=\prod\limits_{z:\widetilde{P}(z)=0}\frac{z^{n}-1}{z-1}= n^{2}\prod\limits_{j=1}^{s-1}\frac{z_{j}^{n}-1}{z_{j}-1}\frac{z_{j}^{-n}-1}{z_{j}^{-1}-1}= n^{2}\prod\limits_{j=1}^{s-1}\frac{z_{j}^{n}+z_{j}^{-n}-2}{z_{j}+z_{j}^{-1}-2}\\
\nonumber &&=n^{2}\prod\limits_{j=1}^{s-1}\frac{2T_{n}(w_{j})-2}{2w_{j}-2} =n^{2}\prod\limits_{j=1}^{s-1}\frac{T_{n}(w_{j})-1}{w_{j}-1}.
\end{eqnarray}
Combine (\ref{tauH}), (\ref{newP}) and (\ref{Hlemma}) we have the following formula for the number of spanning trees \begin{equation}\label{before}
\tau(n)=(-1)^{(m+s)(n-1)}n\,\tau(H)\prod\limits_{j=1}^{s-1}\frac{T_{n}(w_{j})-1}{w_{j}-1}.
\end{equation}

We have the following important statement   from formula (\ref{before}). 

 \bigskip \noindent {\bf Claim:}         \enspace    The number of spanning trees $\tau(n)$ is a multiple of $n\,\tau(H).$

\medskip \noindent{Proof of Claim:}   \hspace{-.5cm}      \enspace  
To prove the lemma we have to show that the number $R=\prod\limits_{j=1}^{s-1}\frac{T_{n}(w_{j})-1}{w_{j}-1}$ is an integer. Indeed, setting $w=\frac{\zeta+2}{2}$ one can represent $R$ in the form 
$$R=\prod\limits_{j=1}^{s-1}\frac{2T_{n}(\frac{\zeta_{j}+2}{2})-2}{\zeta_{j}},$$ 
where $\zeta_{j},\,j=1,2,\ldots,s-1$ are non-zero root of the equation $Q(\frac{\zeta+2}{2})=0.$  We note that  the function $j_n(\zeta)=2T_{n}(\frac{\zeta+2}{2})$ satisfy the recursive relation 
$j_{n+1}(\zeta)=(\zeta+2)j_{n}(\zeta)-j_{n-1}(\zeta)$ 
with initial data $j_0(\zeta)=2$ and $j_1(\zeta)=\zeta+2.$ Hence, $j_{n}(\zeta)$ is a monic polynomial of degree $n$ with integer coefficients. Since $2T_n(1)=2,$ the same is true for the polynomial $f(\zeta)=\frac{2T_{n}(\frac{\zeta+2}{2})-2}{\zeta}.$ By definition, $Q(w)$ is an integer  polynomial in the variables $2T_{s_{i,j}}(w),\,i=1,2,\ldots,m,\,j=1,2,\ldots,k_i.$ By Lemma~\ref{lemma1} we have $Q(1)=0,$ and $Q^{\prime}(1)\neq 0.$ Also, the leading coefficient of $Q(w)$ is equal to $(-1)^{m}2^{s},$ where $s$ is the degree of $Q(w).$ Hence, $g(\zeta)=\frac{1}{\zeta}Q(\frac{\zeta+2}{2})$ with $g(0)\neq 0$ is also a monic polynomial with integer coefficients. Taking this into account, we get  $R={\rm Res}\,(f(\zeta),g(\zeta)).$ Since both $f(\zeta)$ and $g(\zeta)$ are polynomials with integer coefficients, $R$ is integer. {\hfill \qed \newline \smallskip}

Since $\tau(n)$ is a positive number, by (\ref{before}) we obtain

\begin{equation}\label{modulo}
\tau(n)=n\,\tau(H)\prod_{j=1}^{s-1}\left|\frac{T_{n}(w_{j})-1}{w_{j}-1}\right|=
n\,\tau(H)\prod_{j=1}^{s-1}|T_{n}(w_{j})-1|\big/\prod_{j=1}^{s-1}|w_{j}-1|.
\end{equation}

Now we evaluate the product $\prod_{j=1}^{s-1}|w_{j}-1|.$ We note that from Lemma~\ref{lemma1} the polynomial $Q(w)$ has the leading coefficient $a_{0}=(-1)^{m}2^{s},\,Q(1)=0$ and $Q^\prime(1)=-2q,$ where $q=\sum\limits_{i=1}^{m}\sum\limits_{j=1}^{k_{i}}s_{i,j}^2.$

As a result, we have
\begin{eqnarray}\label{eqQ}
\prod_{j=1}^{s-1}|w_{j}-1|=|\frac{1}{a_0}Q^\prime(1)|
=\frac{2q}{2^{s}}=\frac{q}{2^{s-1}}.
\end{eqnarray}

Substituting equation (\ref{eqQ}) into equation (\ref{modulo}) we finish the proof of the theorem.
\end{pf}

\section{Arithmetical properties of complexity for  the graph $H_n$}\label{arithmetic}

Let $H$ be a finite connected graph on $m$ vertices. Consider the circulant foliation $H_n=H_n(G_1,\,G_2,\ldots,G_m),$ where $G_i=C_n(s_{i,1},\,s_{i,2},\ldots,s_{i,k_i}),\,i=1,2,\ldots,m.$ Recall that any positive integer $s$ can be uniquely represented in the form $s=p \,r^2,$ where $p$ and $r$ are positive integers and $p$ is square-free. We will call $p$ the \textit{square-free part} of $s.$

\bigskip

\begin{theorem}\label{lorenzini}
Let $\tau(n)$ be the number of spanning trees in the graph $H_n.$ Denoted by $p$ is the square free parts of   $Q(-1).$ 
Then there exists an integer sequence $a(n)$ such that
\begin{enumerate}
\item[$1^{0}$] $\tau(n)= n\,\tau(H)\,a(n)^{2},$ if $n$ is odd,
\item[$2^{0}$] $\tau(n)=p\,n\,\tau(H)\,a(n)^{2},$ if $n$ is even.
 \end{enumerate}

\end{theorem}

\begin{pf}
By formula~(\ref{tauH}), we have $n\,\tau(n)=\tau(H)\prod_{j=1}^{n-1}\lambda_{1,j}\lambda_{2,j}\ldots\lambda_{m,j}.$ Note that $\lambda_{1,j}\lambda_{2,j}\ldots\lambda_{m,j}=P(\varepsilon_{n}^{j})= P(\varepsilon_{n}^{n-j})=\lambda_{1,n-j}\lambda_{2,n-j}\ldots\lambda_{m,n-j}.$ Define $c(n)=\prod\limits_{j=1}^{\frac{n-1}{2}}\lambda_{1,j}\lambda_{2,j}\ldots\lambda_{m,j},$ if $n$ is odd and $d(n)=\prod\limits_{j=1}^{\frac{n}{2}-1}\lambda_{1,j}\lambda_{2,j}\ldots\lambda_{m,j},$ if $n$ is even. By \cite{Lor}, each algebraic number $\lambda_{i,j}$ comes into the products $\prod_{j=1}^{(n-1)/2}\lambda_{1,j}\lambda_{2,j}\ldots\lambda_{m,j}$ and $\prod_{j=1}^{n/2-1}\lambda_{1,j}\lambda_{2,j}\ldots\lambda_{m,j}$ with all of its Galois conjugate elements. Therefore, both products $c(n)$ and $d(n)$ are integers. Moreover, if $n$ is even we get $\lambda_{1,\frac{n}{2}}\lambda_{2,\frac{n}{2}}\ldots\lambda_{m,\frac{n}{2}}=P(-1)=Q(-1).$ We note that $Q(-1)$  is always a positive integer. The precise   formula for it is given in Remark~1.  

Now, we have $n\tau(n)=\tau(H)\,c(n)^2$ if $n$ is odd, and $n\tau(n)=\tau(H)\,Q(-1)\,d(n)^2$ if $n$ is even. Let $ Q(-1)=p\,r ^{2},$ where $p$ is a  square free number.  Then

\begin{enumerate}
\item[$1^{\circ}$] $\displaystyle{\frac{\tau(n)}{n\,\tau(H)}=\left(\frac{c(n)}{n }\right)^2}$  if $n$    is odd, 
\item[$2^{\circ}$]  $\displaystyle{\frac{\tau(n)}{n\,\tau(H)}=p\left(\frac{r\,d(n)}{n}\right)^2}$ if $n$  is even.  
 
\end{enumerate}

 By Claim in the proof of Theorem~\ref{theorem1},  the quotient $\frac{\tau(n)}{n\,\tau(H)}$  is an integer. 
Since $p$ is square free,  the squared rational numbers in $1^{\circ}$ and $2^{\circ}$ are integer. 
Setting $a(n)=\frac{c(n)}{n}$ in the first case, and $\ a(n)=\frac{r\,d(n)}{n}$ in the second we finish the proof of   the theorem.
\end{pf}

\subsection*{Remark 1}
Denoted by $t_i$ the number of odd elements in the sequence $s_{i,1},\,s_{i,2},\ldots,s_{i,k_i}.$ Then $Q(-1)=\det\,L(H,W),$ where $W=(d_{1}+4t_{1},d_{2}+4t_{2},\ldots,d_{m}+4t_{m}).$ Indeed, $Q(w)=\det\,L(H,W),$ where $W=(w_{1},w_{2},\ldots,w_{m})$ and $w_{i}=2k_{i}+d_{i}-\sum\limits_{j=1}^{k_{i}}2T_{s_{i,j}}(w).$ If $w=-1$ we have $T_{s_{i,j}}(-1)=\cos(s_{i,j}\arccos(-1))=\cos(s_{i,j}\pi)=(-1)^{s_{i,j}}$ and $w_i=d_{i}+4\sum\limits_{j=1}^{k_{i}}\frac{1-(-1)^{s_{i,j}}}{2}=d_{i}+4t_{i}.$

\section{Asymptotic formulas for the number of spanning trees }\label{asymptotic}

In this section we obtain the following result.
\begin{theorem}\label{theorem3}
The asymptotic behaviour for the number of spanning trees $\tau(n)$ in the graph $H_n$ with $\gcd(s_{i,p},\,i=1,\ldots,m,\,p=1,\ldots,k_i)=1$ is given by the formula
$$\tau(n)\sim \frac{n}{q}A^n,\,n\to\infty,$$ where $q=\sum\limits_{i=1}^{m}\sum\limits_{j=1}^{k_i}s_{i,j}^2$ and $A=\exp({\int\limits_{0}^{1}\log|Q(\cos{2 \pi t})|\textrm{d}t}).$
\end{theorem}

To prove the theorem we need the following preliminary lemmas.
\bigskip
\begin{lemma}\label{lemma2}
Let $a_{i,j},(a_{i,i}=0),\,i,j=1,2,\ldots,m$ be non-negative numbers. Let $$D(x_1,x_2,\ldots,x_m)=\det\left(\begin{array}{ccccc}x_1 & -a_{1,2} & -a_{1,3}& \ldots & -a_{1,m} \\
-a_{2,1} & x_2 & -a_{2,3} & \ldots & -a_{2,m} \\
  & \vdots &   & \ddots & \vdots \\
-a_{m,1} & -a_{m,2} & -a_{m,3} & \ldots & x_m\end{array}\right).$$ Then for $x_i\ge d_i=\sum\limits_{j=1}^ma_{i,j},\,i =1,2,\ldots,m$  we have $D(x_1,x_2,\ldots,x_m)\ge0.$  The equality $D(x_1,x_2,\ldots,x_m)=0$ holds if and only if $x_i=d_i,\,i =1,2,\ldots,m.$
\end{lemma}

\begin{pf} We use induction on $m$ to prove the lemma. For $m=1$ we have $D(x_1)=x_1\ge a_{1,1}=0$ and $D(x_1)=0$ iff $x_1=a_{1,1}.$ For $m=2$ one has $D(x_1,x_2)=x_1x_2-a_{1,2}a_{2,1}\ge0$ with $D(x_1,x_2)=0$ if and only if $x_1=a_{1,2}$ and $x_2=a_{2,1}.$  Suppose that $m>2$ and lemma is true for all $D(x_1,x_2,\ldots,x_k)$ with $k<m.$

The $(i,i)$-th minor of the matrix in the statement of lemma is denote by $D(x_1,\ldots,\hat{x}_i,\ldots,x_m),$ where $\hat{x}_i$ means that the variable $x_i$ is dropped. We note that $D^{\prime}_{x_1}(x_1,x_2,\ldots,x_k)=D(x_2,\ldots,x_k).$  Since
$$x_2\ge d_2=\sum\limits_{j=1}^ma_{2,j}\ge\sum\limits_{j=2}^ma_{2,j},x_3\ge d_3=\sum\limits_{j=1}^ma_{3,j}\ge\sum\limits_{j=2}^ma_{3,j}, \ldots,x_m\ge d_m=\sum\limits_{j=1}^ma_{m,j}\ge\sum\limits_{j=2}^ma_{m,j},$$ the function $D(x_2,\ldots,x_m)$ satisfies the conditions of lemma. Hence, $D(x_2,\ldots,x_m)\ge 0.$ In a similar way, for $i=2,\ldots,m$ we have
$$D^{\prime}_{x_i}(x_1,x_2,\ldots,x_m)=D(x_1,\ldots,\hat{x}_i\ldots,x_m)\ge0.$$ Since $D(d_1,d_2,\ldots,d_m)=0,$ we obtain $D(x_1,x_2,\ldots,x_m)\ge0$ for all $x_i\ge d_i,\,i =1,2,\ldots,m.$ If for some $i_0$ we have $x_{i_0}>d_{i_0},$ then, by induction, for all $i\neq i_{0}$ we get $D^{\prime}_{x_{i_0}}(x_1,x_2,\ldots,x_m)=D(x_2,\ldots,\hat{x}_{i_0}\ldots,x_m)>0$ and $D(x_1,x_2,\ldots,x_m)>0.$
\end{pf}

\begin{lemma}\label{nounityroots}
Let $\gcd(s_{i,j},\,i=1,\ldots,m,\,j=1,\ldots,k_i)=1$ and $s=s_{1,k_1}+s_{2,k_2}\ldots+s_{m,k_m}$ Then the roots of the Laurent polynomial $P(z)$ counted with multiplicities are $1,\,1,\,z_{1},\,1/z_{1},\ldots,\,z_{s-1},\,1/z_{s-1},$ where we have $|z_{p}|\neq1,\,p=1,2,\ldots,s-1.$ Polynomial $Q(w)$ has the roots $1,\,w_{1},\ldots, w_{s-1},$ where $w_{p}=\frac{1}{2}(z_{p}+z_{p}^{-1})$ for all $p=1,\,2,\ldots,s-1.$
\end{lemma}

\begin{pf} By Lemma~\ref{lemma1} we have $P(z)=Q(\frac{1}{2}(z+z^{-1}))$ and $Q(w)$ has the simple root $w=1.$

Since the mapping $w=\frac{1}{2}(z+z^{-1})$ is two-to-one, the Laurent polynomial $P(z)$ has the double root $z=1.$

To prove the lemma we suppose that the Laurent polynomial $P(z)$ has a root $z_{0}$ such that $|z_{0}|=1$ and $z_{0}\neq1.$ Then $z_{0}=e^{\textrm{i}\,\varphi_{0}},\,\varphi_{0}\in\mathbb{R}\setminus2\pi\mathbb{Z}.$ Now we have
$$P(e^{\textrm{i}\,\varphi_{0}})=
\det\left(\begin{array}{ccccc}x_1 & -a_{1,2} & -a_{1,3}& \ldots & -a_{1,m} \\
-a_{2,1} & x_2 & -a_{2,3} & \ldots & -a_{2,m} \\
  & \vdots &   & \ddots & \vdots \\
-a_{m,1} & -a_{m,2} & -a_{m,3} & \ldots & x_m\end{array}\right),$$ where $$x_{i}=2k_{i}+d_{i}-\sum\limits_{j=1}^{k_{i}}(z_{0}^{s_{i,j}}+z_{0}^{-s_{i,j}}) =2k_{i}+d_{i}-\sum\limits_{j=1}^{k_{i}}(e^{\textrm{i}\,s_{i,j}\,\varphi_{0}}+e^{-\textrm{i}\,s_{i,j}\,\varphi_{0}}) =d_{i}+\sum\limits_{j=1}^{k_{i}}(2-2\cos(s_{i,j}\,\varphi_{0})).$$

Since $d_{i}=\sum\limits_{j=1}^{m}a_{i,j}$ and $x_{i}\geq d_{i},$ the conditions of Lemma~\ref{lemma2} are satisfied. Hence $P(e^{\textrm{i}\,\varphi_{0}})=0$ if and only if $x_{i}=d_{i},\,i=1,\ldots,m.$ Then $\cos(s_{i,j}\,\varphi_{0})=1$ for all $i=1,\ldots,m,\,j=1,\ldots,k_{i}.$ So $s_{i,j}\,\varphi_{0}=2\pi m_{i,j}$ for some integer $m_{i,j}.$ As $\gcd(s_{i,j},\,i=1,\ldots,m,\,j=1,\ldots,k_i)=1$ there exist integers $p_{i,j}$ such that $\sum\limits_{i=1}^{m}\sum\limits_{j=1}^{k_{i}}s_{i,j}p_{i,j}=1.$ See, for example, (\cite{Apost}, p.~21). Hence, $\varphi_0=\varphi_0\sum\limits_{i=1}^{m}\sum\limits_{j=1}^{k_{i}}s_{i,j}p_{i,j}= 2\pi\sum\limits_{i=1}^{m}\sum\limits_{j=1}^{k_{i}}m_{i,j}p_{i,j}\in2\pi\mathbb{Z}.$ Contradiction.
\end{pf}

Now we come to the proof of the Theorem~\ref{theorem3}

\begin{pf}
By theorem \ref{theorem1} we have $\tau(n)=\frac{n\tau(H)}{q}\prod\limits_{j=1}^{s-1}|{2T_{n}(w_{j})-2}|,$ where $q=\sum\limits_{i=1}^{m}\sum\limits_{j=1}^{k_{i}}s_{i,j}^2$ and $w_{j},\,j=1,2,\ldots,s-1$ are roots of the polynomial $Q(w)$ different from $1$.

By lemma~\ref{nounityroots}, $T_{n}(w_{j})=\frac{1}{2}(z_{j}^{n}+z_{j}^{-n}),$ where the $z_{j}$ and $1/z_{j}$ are roots of the polynomial $P(z)$ with the property $|z_{j}|\neq1,\,j=1,2,\ldots,s-1.$ Replacing $z_{j}$ by $1/z_{j},$ if it is necessary, we can assume that  $|z_j|>1$ for all $j=1,2,\ldots,s-1.$ Then $T_{n}(w_{j})\sim\frac{1}{2}z_{j}^{n}$ and $|2T_{n}(w_{s})-2|\sim|z_{s}|^{n}$ as $n\to\infty.$ Hence
$$\frac{n\tau(H)}{q}\prod_{j=1}^{s-1}|2T_{n}(w_{j})-2|\sim\frac{n\tau(H)}{q}\prod_{j=1}^{s-1} |z_{j}|^{n}=\frac{n\tau(H)}{q}\prod\limits_{P(z)=0,\,|z|>1}|z|^{n}=\frac{n A^n}{q},$$ where $A=\prod\limits_{P(z)=0,\,|z|>1}|z|$ is the Mahler measure of the polynomial $P(z).$ By (\cite{EverWard}, p.~67), we have $A=\exp\left(\int_{0}^{1}\log|P(e^{2 \pi i t })|\textrm{d}t\right).$ Since $P(z)=Q(\frac{1}{2}(z+z^{-1})),$ we get $A=\exp({\int\limits_{0}^{1}\log|Q(\cos{2 \pi t})|\textrm{d}t}).$ 
The theorem is proved.
\end{pf}

\subsection*{Remark 2}

We note that $Q(\cos(2\pi t))=\det\,L(H,W),$ where $W=(w_{1},w_{2},\ldots,w_{m})$ and $w_{i}=2k_{i}+d_{i}-\sum\limits_{j=1}^{k_{i}}2T_{s_{i,j}}(\cos(2\pi t))=d_{i}+4\sum\limits_{j=1}^{k_{i}}\sin^2(s_{i,j}\pi t),\,i=1,2,\ldots,m.$

\section{Examples}\label{examples}

\subsection{Circulant graph $C_{n}(s_{1},s_{2},\ldots,s_{k}).$}\label{example1}

We consider the classical circulant graph $C_{n}(s_{1},s_{2},\ldots,s_{k})$ as a foliation $H_{n}(G_{1})$
on the one vertex graph $H=\{v_{1}\}$ with the fiber $G_{1}=C_{n}(s_{1},s_{2},\ldots,s_{k}).$
In this case $d_{1}=0,\,L(H,X)=(x_{1}),\,P(z)=2k-\sum\limits_{p=1}^{k}(z^{s_p}+z^{-s_p})$ and
its Chebyshev transform is $Q(w)=2k-\sum\limits_{p=1}^{k}2T_{s_p}(w).$
Different aspects of complexity for circulant graphs were investigated in the papers \cite{XiebinLinZhang, ZhangYongGolin, Golin2010, MedMed18, MedMed17}.

\subsection{$I$-graph $I(n,k,l)$ and the generalized Petersen graph $GP(n,k).$}

Let $H$ be a path graph on two vertices, $G_{1}=C_{n}(k)$ and $G_{2}=C_{n}(l).$ Then $I(n,k,l)=H_{n}(G_{1},G_{2})$ and $GP(n,k)=I(n,k,1).$ We get $P(z)=(3-z^{k}-z^{-k})(3-z^{l}-z^{-l})-1$ and $Q(w)=(3-2T_{k}(w))(3-2T_{l}(w))-1.$ The arithmetical and asymptotical properties of complexity for $I$-graphs were studied in \cite{Ilya}.

\subsection{Sandwich of $m$ circulant graphs.}

Consider a path graph $H$ on $m$ vertices.
Then $H_{n}(G_{1},G_{2},\ldots,G_{m})$ is a {\it sandwich graph} of circulant graphs
$G_{1},G_{2},\ldots,G_{m}.$ Here  $d_{1}=d_{m}=1$ and $d_{i}=2,\,i=2,\ldots,m-1.$
We set $$D(x_1,x_2,\ldots,x_m)=
\det\left(\begin{array}{ccccccc}
x_1 & -1  & 0  & \ldots & 0 & 0& 0 \\
-1  & x_2 & -1 & \ldots & 0 & 0& 0 \\
& \vdots & & \ddots& & \vdots & \\
0 & 0 & 0 & \ldots & -1 & x_{m-1} & -1\\
0 & 0 & 0 & \ldots &  0 & -1 & x_m\end{array}\right).$$   By direct calculation we obtain
$$D(x_1,x_2,\ldots,x_m)=x_1D(x_2,\ldots,x_m)-D(x_3,\ldots,x_m),\, D(x_1)=x_1,\,D(x_1,x_2)=x_1x_2-1.$$ Then
$Q(w)=D(w_1,w_2,\ldots,w_m)$  and  $Q(-1)=D(d_1+4t_1,d_2+4t_2,\ldots, d_m+4t_m),$
where $w_i$ and $t_i$ are the same as in Theorem~\ref{theorem3}. 

\subsection{Generalized $Y$-graph.}

Consider the generalized $Y$-graph $Y_{n}(G_{1},G_{2},G_{3})$ where $G_{i}=C_{n}(s_{i,1},\,s_{i,2},\ldots,s_{i,k_i}),\,i=1,2,3.$ Here 
$$Q(w)=3A_{1}(w)A_{2}(w)A_{3}(w)-A_{1}(w)A_{2}(w)-A_{1}(w)A_{3}(w)-A_{2}(w)A_{3}(w),$$ 
where $A_{i}(w)=2k_{i}+1-\sum\limits_{j=1}^{k_{i}}2T_{s_{i,j}}(w).$ 

\subsection{Generalized $H$-graph.}

Consider the generalized $H$-graph $H_{n}(G_{1},G_{2},G_{3},G_{4}),$ where $G_{i}=C_{n}(s_{i,1},\,s_{i,2},\ldots,s_{i,k_i}),\,i=1,2,3,4.$ Now we have 
$$Q(w)=A_{1}(w)A_{2}(w)A_{3}(w)A_{4}(w)\left((3-\frac{1}{A_{1}(w)}-\frac{1}{A_{2}(w)})
(3-\frac{1}{A_{3}(w)}-\frac{1}{A_{4}(w)})-1\right)$$ where $A_{i}(w)$ are the same as above.

\subsection{Discrete torus $T_{n,m}=C_n\times C_m.$}

We have $T_{n,m}=H_{n}(\underbrace{C_{n}(1),\ldots,C_{n}(1)}_{m \textrm{ times}}),$ where $H=C_{m}(1)$ is the cyclic graph on $n$ vertices. So, the generalized Laplacian matrix  with respect to the set of variables $X=(\underbrace{x,\ldots,x}_{m \textrm{ times}})$ has the form
$L(H,X)=
\left(\begin{array}{cccccc}
x & -1 & 0 & \ldots & 0& -1 \\
-1 & x & -1 & \ldots & 0& 0 \\
& \vdots &  & \ddots & &\vdots \\
-1 & 0 & 0 & \ldots & -1 & x\\
\end{array}\right).$ Then
$L(H,X)$ is an $m\times m$ circulant matrix with eigenvalues $\mu_j=x-e^{\frac{2\pi i j}{m}}-(e^{\frac{2\pi i j}{m}})^{m-1}=x-2\cos(\frac{2\pi j}{m}),j=0,\ldots,m-1.$ Hence, $\det L(H,X)=\prod\limits_{j=0}^{m-1}\mu_j=2 T_{m}(x/2)-2.$ Substituting $x=4-z-z^{-1}$ and $w=\frac{1}{2}(z+z^{-1}),$ we get $Q(w)=2 T_{m}(2-w)-2.$

\subsection{Direct product $C_n\times H$ where $H$ is a regular graph.}

Let $H$ be a connected $d$-regular graph. One can identify the direct product $C_n\times H$ with $H_{n}=H_{n}(\underbrace{C_{n}(1),\ldots,C_{n}(1)}_{m \textrm{ times}}).$ Let $X=(\underbrace{x,\ldots,x}_{m \textrm{ times}})$. Now $L(H,X)=x I_{m}-A(H).$ Hence, $\det\,L(H,X)$ coincides with the characteristic polynomial $\chi_{H}(x)$ of graph $H.$ We have $Q(w)=\chi_{H}(2+d-2w).$ Then $Q(-1)=\chi_{H}(4+d).$

\section*{ACKNOWLEDGMENTS}
The work was partially supported by the Korean-Russian bilateral project. The first author was supported in part by the Basic Science Research Program through the National Research Foundation of Korea (NRF) funded by the Ministry of Education (2018R1D1A1B05048450).  The second and the third authors
were partially supported by the Russian Foundation for Basic Research (grants 18-01-00420 and 18-501-51021).  The results given  in Sections 5 and 6 are supported by the Laboratory of Topology and Dynamics, Novosibirsk State University  (contract no. 14.Y26.31.0025 with the Ministry of Education and Science of the Russian Federation).

\end{document}